\documentclass{amsart}

\newtheorem{prop}{Proposition}
\newtheorem{theo}{Theorem}
\newtheorem{Lemma}{Lemma}
\newtheorem{cor}{Corollary}

\def\Bbb{\mathbb}

\newcommand{\na}{\nabla}

\newcommand{\om}{\omega}
\newcommand{\Om}{\Omega}
\newcommand{\la}{\lambda}
\newcommand{\La}{\Lambda}

\newcommand{\ka}{K{\"a}hler }

\newcommand{\ranti}{{\rm Ric}''}

\newcommand{\rhs}{\rho^{*}}

\newcommand{\rhssym}{(\rho^{*})'}
\newcommand{\rhsskew}{(\rho^{*})''}

\newcommand{\leftr}{[\hbox{\hspace{-0.15em}}[}
\newcommand{\rightr}{]\hbox{\hspace{-0.15em}}]}
\newcommand{\ocap}{\bigcirc \hspace{-10.5pt} \wedge\hspace{0.1cm}}
\newcommand{\ocaptext}{\bigcirc \hspace{-8.5pt} \wedge\hspace{0.1cm}}

\title[K{\"a}hler manifolds with Ricci tensor of constant eigenvalues]{A
splitting theorem for K{\"a}hler  manifolds whose Ricci tensors
have constant eigenvalues}
\author{Vestislav Apostolov, Tedi Dr\u{a}ghici and Andrei Moroianu}
\thanks{The first two authors were supported in part by
NSF grant INT-9903302}
\address{Vestislav Apostolov\\ D\'epartement de Math\'ematiques \\
Universit\'e du Qu\'ebec \`a Montr\'eal \\
Case postale 8888 \\ succursale centre-ville\\ Montr\'eal (Qu\'ebec)\\
H3C 3P8, Canada}
\email{apostolo@math.uqam.ca}

\address{Tedi Dr\u{a}ghici \\ Department of Mathematics \\ Florida
International University \\ Miami FL 33199 \\ USA}
\email{draghici@fiu.edu}

\address{Andrei Moroianu \\ CMAT\\ {\'E}cole Polytechnique \\ UMR 7640 du CNRS
\\ 91128 Palaiseau \\ France}
\email{am@math.polytechnique.fr}
\begin{document}

\begin{abstract}
It is proved that a compact K{\"a}hler manifold whose Ricci
tensor has two distinct constant non-negative eigenvalues is
locally the product of two K{\"a}hler-Einstein manifolds. A
stronger result is established for the case of K\"ahler surfaces.
Without the compactness assumption, irreducible K\"ahler manifolds
with Ricci tensor having two distinct constant eigenvalues are
shown to exist in various situations: there are homogeneous
examples of any complex dimension $n\ge 2$ with one eigenvalue
negative and the other one positive or zero; there are
homogeneous examples of any complex dimension $n\ge 3$ with two
negative eigenvalues; there are non-homogeneous examples of
complex dimension 2 with one of the eigenvalues zero. The problem
of existence of K{\"a}hler metrics whose Ricci tensor has two
distinct constant eigenvalues is related to the celebrated (still
open) conjecture of Goldberg \cite{goldberg}. Consequently, the
irreducible homogeneous examples with negative eigenvalues give
rise to {\it complete} Einstein strictly almost K\"ahler metrics
of any even real dimension greater than 4.

\vspace{0.1cm}
\end{abstract}

\maketitle

\section{Introduction}

In this note we  consider compact
K{\"a}hler manifolds $(M,g,J)$ whose Ricci tensor, ${\rm Ric}$,
(or rather the symmetric endomorphism corresponding to ${\rm Ric}$
via the metric)
has two distinct constant {eigenvalues}.
Since for K{\"a}hler manifolds  the Ricci tensor is
invariant under the action of $J$ (i.e. it satisfies
${\rm Ric}(J\cdot, J\cdot)= {\rm Ric}(\cdot, \cdot)$),
each eigenvalue is of even multiplicity.

One source of motivation for studying such manifolds comes from a
conjecture of Lichnerowicz concerning
compact K{\"a}hler spin manifolds with least possible (compared to
the scalar curvature) eigenvalue of the Dirac operator. It was shown
by Kirchberg in 1990 that every eigenvalue $\nu$ of the Dirac
operator on a compact K{\"a}hler spin manifold $M$ of even complex
dimension $n$ satisfies
$$\nu^2\ge \frac{n}{4(n-1)}\inf_M {s},$$
where $s$ denotes the scalar curvature of $M$. The limiting case of
this inequality is characterized by the existence of spinor fields
satisfying a certain first order differential equation, and it was
conjectured by Lichnerowicz that the Ricci tensor of such manifolds
has to be parallel. In 1997, it was shown
in \cite{And} that the Ricci tensor of every limiting manifold for the
above inequality has  two constant eigenvalues, one positive with
multiplicity $2n-2$ and one zero with multiplicity 2;
after some more work, the conjecture was proved in \cite{andrej}
using Spin$^c$ geometry.

Nevertheless, the general question whether the Ricci tensor of a compact
K{\"a}hler manifold has to be parallel as soon as it has two
non-negative constant eigenvalues was left open in \cite{andrej}.
A positive
answer to this question is summarized in the following
splitting theorem.

\begin{theo}\label{th0} Let $(M,g,J)$ be a compact
K{\"a}hler manifold whose Ricci tensor has two distinct constant non-negative
eigenvalues $\la$ and $\mu$.
Then the universal cover of $(M,g,J)$ is the product of two
simply connected K{\"a}hler--Einstein manifolds of
scalar curvatures $\la$ and $\mu$, respectively.
\end{theo}

Note that {\it local} irreducible examples of
K{\"a}hler manifolds  with eigenvalues of the Ricci tensor
equal to 0 and 1 are known to exist in complex dimension two,
cf. \cite{Br} and \cite[Remark 1(c)]{A-A-D}. Note also that the above
result fails if one allows the Ricci tensor to have more than two
different eigenvalues, as shown by the compact homogeneous K{\"a}hler
manifolds (the generalized complex flag manifolds).

We also mention that K\"ahler manifolds with Ricci tensors having constant eigenvalues
are related to different generalizations of
Calabi's construction of K\"ahler metrics with  constant scalar
curvature \cite{hwang-singer}; in this vein R. Bielawski recently
proved \cite{bielawski} that the total space of the canonical
bundle of a K\"ahler manifold with Ricci tensor having constant
eigenvalues carries a Ricci-flat K\"ahler metric. Interesting
new  examples actually appear when the underlying manifold is
irreducible and non-homogeneous; the {\it completeness} of the
Ricci-flat metric corresponds to the case when the underlying
manifold is {\it compact} and the Ricci tensor has {\it positive}
constant eigenvalues.

Another motivation for our study
came from an a priori unexpected link with a conjecture of
Goldberg \cite{goldberg}, which states that any compact {\it
Einstein} almost K{\"a}hler manifold is, in fact,
K{\"a}hler--Einstein. This link is presented in the first part of
the Section 2, while preparing the ground for the proof of
Theorem \ref{th0}. For any K{\"a}hler manifold $(M,g,J)$ with
Ricci tensor having two distinct constant eigenvalues, one can
define another $g$-orthogonal almost complex structure ${\bar
J}$, by changing the sign of $J$ on one of the eigenspaces of
${\rm Ric}$. The new almost complex structure ${\bar J}$, which
is not integrable in general, commutes with $J$ and has a closed
fundamental 2-form, i.e.  $(g,{\bar J})$ gives rise to an almost
K{\"a}hler structure on the manifold. The integrability of ${\bar
J}$ holds precisely when the Ricci tensor of $g$ is parallel, or
equivalently, when $g$ is locally a product of two
K{\"a}hler--Einstein metrics; see Lemma \ref{lem1}. Moreover, any
K{\"a}hler structure $(g,J)$ with Ricci tensor having two
distinct constant eigenvalues either both positive, or both
negative, determines (and is determined by) a certain Einstein
almost K{\"a}hler  structure $({\widetilde g}, {\bar J})$, see
Corollary \ref{cor1}.

\vspace{0.3cm} The proof of Theorem \ref{th0} is completed in
Section 2. We show the integrability of the almost K{\"a}hler
structure $(g, {\bar J})$ using essentially an integral formula
found by Sekigawa in \cite{sek}, where he gave an affirmative
answer to the Goldberg conjecture in the case of non-negative
scalar curvature. However, in our situation a more detailed
analysis  is required to extract the conclusion, which we
accomplish by making use of some Weitzenb{\"o}ck formulae
elaborated in \cite{Bour}.

\vspace{0.3cm}

In Section 3 we consider the {\it homogeneous} K{\"a}hler
manifolds, as an important class of K{\"a}hler manifolds whose
Ricci tensors have constant eigenvalues. Theorem \ref{th0} implies
the splitting of any compact homogeneous K{\"a}hler manifold with
Ricci tensor having two distinct positive eigenvalues (see
Corollary \ref{cor2}); the splitting in the case when one of the
eigenvalues is zero and the other one is positive was already
known from classical results \cite{cheeger-gromoll,Nak}.

On the other hand, according to the structure theorem for
homogeneous \ka manifolds \cite{V-G-PSh,Dor-Nak}, there are lots
of irreducible examples with Ricci tensor having eigenvalues
$(-1,+1)$, or $(-1,0)$.

We finally exhibit irreducible (non-compact) homogeneous
K{\"a}hler manifolds of any complex dimension $n\ge3$ whose Ricci
tensors have two
distinct eigenvalues, both negative and constant. Our
construction relies on an appropriate realization of a bounded
homogeneous domain as a Siegel domain of type {\rm II}; the
metric is then invariant under the simply transitive action of a
solvable group of affine transformations of the relevant Siegel
domain. By virtue of Corollary \ref{cor1}, any such homogeneous
metric can be deformed to a K{\"a}hler--Einstein one, which also
admits a {\it strictly} almost K{\"a}hler structure. (Here and
henceforth, {\it strictly} almost K{\"a}hler structure means that
the corresponding almost complex structure is not integrable). We
thus provide irreducible {\it complete} Einstein strictly almost
K{\"a}hler manifolds of {\it any} even dimension greater than
four. (The previously known such examples were constructed by D.
Alekseevsky  on certain solvable Lie groups of high dimensions,
see {\it e.g.} \cite[pp. 419--421 \& Rem. 14.100]{besse}). Even
more surprising is the fact that some non-compact Hermitian
symmetric spaces admit a strictly almost K{\"a}hler structure commuting with the standard
Hermitian structure
(Corollary \ref{examples}). Particular examples (filling all real
dimensions $2n, n\ge 3$) are the Hermitian symmetric spaces
$M^{2n} = {\rm SO}(2,n)/({\rm SO}(2)\times {\rm SO}(n))$, $n\ge
3$; see Example 1. By contrast, J. Armstrong \cite{Arm} showed that the (unique)
irreducible non-compact Hermitian symmetric space of complex
dimension 2, the complex hyperbolic space, does not admit (even
locally) strictly almost K{\"a}hler structures commuting with the
standard K{\"a}hler structure.
Note that symmetric spaces of {\it
compact type} of any dimension do not (locally) carry  orthogonal strictly
almost K{\"a}hler structures, cf. \cite{Ga}, \cite{her-lamoneda}
(see also Remark 2 below); it is also known from
\cite{Arm,Og-Sek,Olszak} that real hyperbolic spaces do not (locally) admit
orthogonal strictly almost K{\"a}hler structures either.

\vspace{0.3cm}
In the last part of the paper, we provide a stronger version of Theorem
\ref{th0}
for the case of compact K{\"a}hler manifolds of complex dimension 2.

\begin{theo}\label{th1} Let $(M,g,J)$ be a compact K{\"a}hler surface
whose Ricci tensor has two distinct constant eigenvalues. Then one of the
following alternatives holds~:
\begin{enumerate}
\item[{\rm (i)}]  $(M,g,J)$ is locally symmetric, i.e. locally is the
product of Riemann surfaces of distinct constant Gauss curvatures;
\item[{\rm (ii)}] if $(M,g,J)$ is not as described in (i), then
the eigenvalues of the Ricci tensor are both negative and
$(M,J)$ must be a minimal surface of general type with ample
canonical bundle and with even and positive signature. Moreover, in this
case, reversing the orientation, the manifold would admit an Einstein,
strictly almost K{\"a}hler metric.
\end{enumerate}
\end{theo}

The additional ingredients for the proof of Theorem \ref{th1}
come from the classification results of K\"ahler surfaces and
from consequences of the Seiberg-Witten theory, see \cite{Kots,leung,Pe}.
At this time, we do not know if the alternative (ii) may really hold.
As pointed out, if it does, it also provides a counter-example to
the 4-dimensional Goldberg conjecture.

Let us also note here that the compactness assumption in Theorem
\ref{th1} is essential: Complete (non-compact) examples are the
non-symmetric homogeneous \ka surfaces  with Ricci tensor having
eigenvalues $(-1,+1)$ or $(-1,0)$, see respectively \cite{Shima1}
and \cite{Kow}. The case of eigenvalues  $(-1,0)$ corresponds to
the unique proper 3-symmetric space in four dimensions, also
viewed as the non-compact K{\"a}hler geometry  ${\bf F}^4$ of \cite{wall}.
Moreover, the constructions of \cite{Br} and \cite{A-A-D} give
further local examples of non-homogeneous  K{\"a}hler surfaces
with Ricci tensor of constant eigenvalues; in all these examples
one of the two eigenvalues of ${\rm Ric}$ is non-negative,
therefore,  by Theorem \ref{th1},  none of them can be
compactified.

\section{Proof of Theorem \ref{th0}}
\subsection{The commuting almost K{\"a}hler structure}
The main idea of the proof of our results is to construct  an
almost K{\"a}hler structure ${\bar J}$ on ${(M,g,J)}$, which is compatible
with $g$ and commutes with $J$.

\begin{Lemma}\label{lem1} Let $(M,g,J)$ be a  K{\"a}hler manifold whose Ricci
tensor has constant
eigenvalues $\la < \mu$. Denote by $E_{\la}$ and $E_{\mu}$ the corresponding
$J$-invariant
eigenspaces and define a $g$-orthogonal almost complex structure
${\bar J}$ by
setting
${\bar J}_{| E_{\la}}=  J_{| E_{\la}}; \ {\bar J}_{| E_{\mu}}= -J_{| E_{\mu}}$.
Then $J$ and ${\bar J}$ mutually commute and
$(g, {\bar J})$ is an almost K{\"a}hler structure, i.e. the fundamental form
${\overline \Om}(\cdot, \cdot)=g({\bar J}\cdot,\cdot)$
is symplectic.   Moreover, $(g,{\bar J})$  is
K{\"a}hler {\rm (}i.e. ${\bar J}$ is integrable{\rm )} if and only if
$(M,g)$ is locally product of two K{\"a}hler--Einstein manifolds of scalar
curvatures $\la$ and $\mu$, respectively.
\end{Lemma}

\noindent
{\it Proof.}  Denote by $\Omega(\cdot,\cdot)=g(J\cdot,\cdot)$ the
fundamental form of $(g,J)$
and consider the (1,1)-forms $\alpha$ and $\beta$ defined by
$$\alpha (X,Y) = \Om({\rm pr}^{\la}(X),{\rm pr}^{\la}(Y)), \
\forall X,Y \in TM; $$
$$\beta = \Om - \alpha, $$
where ${\rm pr}^{\la}$ (resp. ${\rm pr}^{\mu}$) denotes the orthogonal
projection of the tangent bundle $TM$
onto $E_{\la}$ (resp. $E_{\mu}$).
The Ricci form  $\rho(\cdot,\cdot)= {\rm Ric}(J\cdot, \cdot)$ of $(M,g,J)$ is
then given by
$$\rho = \la \alpha + \mu \beta. $$
As $\Om = \alpha + \beta$ and $\rho$ are both closed (1,1)-forms,
so are the 2-forms $\alpha$ and $\beta$.
By the very definition of ${\bar J}$,
the fundamental form ${\overline \Omega}(\cdot, \cdot)= g({\bar J}\cdot,
\cdot)$
is given by
$${\overline \Omega} = \alpha - \beta, $$
and hence is closed, i.e. $(g,{\bar J},{\overline \Omega})$ is an
almost K{\"a}hler structure; it is \ka  as soon as the Ricci
tensor is parallel (equivalently, $\alpha$ and $\beta$ are
parallel), i.e. when $(M,g)$ is locally a product of two
K{\"a}hler-Einstein manifolds with scalar curvatures $\la$ and
$\mu$, respectively. {\bf Q.E.D.}

\vspace{0.2cm}
\noindent
{\it Remark 1.}  Note that the Ricci tensor of
the almost K{\"a}hler structure
$(g,{\bar J},{\overline \Om})$ constructed in Lemma \ref{lem1} is
${\bar J}$-invariant, i.e. is of type $(1,1)$ with respect to ${\bar J}$.
In the compact case,
almost \ka structures  with Ricci tensors
of type $(1,1)$ are the {critical} points
of the
Hilbert functional, the integral of the scalar curvature, restricted
to the set of all compatible metrics to  symplectic form ${\overline \Om}$ (cf. \cite{BI}).
Compatible
\ka metrics provide absolute maxima for the functional in this
setting and it was a natural question~\cite{BI} to ask if every critical
metric is necessarily K{\"a}hler.  The answer turns out to be
negative in dimension greater than four~\cite{DM}, while in four dimensions
the problem is still open, as no examples of
{\it compact}, non-K{\"a}hler, almost \ka structures with  Ricci
tensor of type $(1,1)$ are known yet.

\vspace{0.2cm}

Let us consider
for a moment the more general context of K{\"a}hler manifolds $(M,g,J)$
which admit a {\it commuting} almost K{\"a}hler structure ${\bar J}$.
Any $g$-orthogonal almost complex
structure ${\bar J}$ which  commutes with and differs from $\pm J$ gives
rise to a
$g$-orthogonal, $J$-invariant endomorphism $Q= - J\circ {\bar J}$ of $TM$
such that $Q^2={\rm Id}_{|_{TM}}$; we thus
define an orthogonal, $J$-invariant splitting of the tangent bundle
$TM$
$$TM= E_{+}\oplus E_{-}$$
into the sum of the $\pm 1$-eigenspaces of $Q$, the (complex) sub-bundles
$E_{\pm}$, respectively. As in the proof of Lemma \ref{lem1}, we consider the
(1,1)-forms $\alpha$ and $\beta$, the restrictions of the fundamental form
$\Om$
of $(g,J)$ to the spaces $E_{+}$ and $E_{-}$, respectively.
The fundamental
forms $\Om$ and ${\overline \Om}$ of $(g,J)$ and $(g,{\bar J})$ are then
given by
$$\Om = \alpha + \beta; \ \ {\overline \Om}= \alpha - \beta,$$
proving that  $\alpha$ and $\beta$ are {\it closed}.
Therefore, corresponding to any
K{\"a}hler metric $(g,J)$ admitting a commuting almost K{\"a}hler structure
${\bar J}$, we may consider a natural 1-parameter family $g^t$
of metrics having the same property  (see \cite{andrej}):
\begin{equation}\label{gt}
g^t=g_{|_{E_{+}}} + t g_{|_{E_{-}}},  \ \ t>0,
\end{equation}
where $g_{|_{E_{+}}}$ {\rm (}resp. $g_{|_{E_{-}}}${\rm })
denotes the restriction of
$g$ to the eigenspaces $E_{+}$ {\rm (}resp. to $E_{-}${\rm )}.
\begin{Lemma}\label{lem2}
For any $t>0$, the metric $g^t$ is
K{\"a}hler with respect to $J$, almost K{\"a}hler with respect to ${\bar J}$
and has the same Ricci tensor as the metric $g = g^1$.
\end{Lemma}
\noindent {\it Proof.} The first statements follow from the fact
that the fundamental form of $(g^t, J)$ (resp. $(g^t,{\bar J})$)
is closed as being equal to $\alpha + t \beta$ (resp. $\alpha - t
\beta$), where $\alpha$ and $\beta$ are constructed as above with
respect to $g=g^1$. For the last claim, note that the volume form
of the metric $g^t$ is a constant multiple of the volume form of
$g=g^1$, so, from the local expression in complex coordinates,
the Ricci forms of the K{\"a}hler structures $(g^t,J)$ and $(g,J)$
coincide. {\bf Q.E.D.}

\vspace{0.2cm} \noindent As for K{\"a}hler metrics with Ricci
tensor having distinct constant eigenvalues, Lemma \ref{lem2}
shows that one can deform any such given metric to one whose
Ricci tensor has constant eigenvalues equal to $-1,0$ or $+1$. In
particular, we get

\begin{cor} \label{cor1} On a complex manifold $(M^{2n},J)$ there is a one-to-one
correspondence between K{\"a}hler metrics with Ricci tensor of
constant eigenvalues $\la < \mu$ with $\la\mu >0$ and
K{\"a}hler-Einstein metrics ${\widetilde g}$ of scalar curvature
$2n\la$ carrying an orthogonal almost K{\"a}hler structure  ${\bar
J}$ which commutes with and differs from $\pm J$; in this
correspondence ${\bar J}$ is  compatible also with $g$ and
coincides {\rm (}up to sign{\rm )} with the almost K{\"a}hler
structure  defined in Lemma \ref{lem1};   moreover, ${\bar J}$ is
integrable precisely when $g$ {\rm (}and ${\widetilde g}${\rm )}
is locally product of two K{\"a}hler--Einstein metrics.
\end{cor}
\noindent {\it Proof.} Let $(M,g,J)$ be a K{\"a}hler manifold
whose Ricci tensor has constant eigenvalues $\la < \mu$ and
${\bar J}$ be the almost K{\"a}hler structure commuting with $J$
given by Lemma \ref{lem1}. It is easy to see that the $\pm
1$-eigenspaces of the endomorphism $Q= - J\circ {\bar J}$ above
are given by $E_+ = E_{\la}, \ E_-=E_{\mu}$, where, we recall,
$E_{\la}$ and $E_{\mu}$ are the eigenspaces of ${\rm Ric}$. By
Lemma 2, the metric ${\widetilde g} = g^{\mu /\la}$ obtained via
(\ref{gt}) is K{\"a}hler--Einstein with scalar curvature $2n\la$.
Conversely, starting from a K{\"a}hler--Einstein structure
$({\widetilde g}, J)$ of scalar curvature $2n\la$, endowed with an
almost K{\"a}hler structure ${\bar J}$ commuting with $J$, the
deformation (\ref{gt}) provides a K{\"a}hler metric $(g,J)$ whose
Ricci tensor has constant eigenvalues $\la <\mu$, by putting
$g^1={\widetilde g}$ and $g= g^{\la / \mu}$. The almost complex
structure ${\bar J}$ is compatible to both $g$ and ${\widetilde
g}$. It is clear then that the common Ricci tensor of $g$ and
${\widetilde g}$ is ${\bar J}$-invariant, and therefore, ${\bar
J}$ coincides (up to sign) with the almost complex structure
defined in Lemma \ref{lem1}. By Lemma \ref{lem1} we also conclude
that the integrability of ${\bar J}$ is equivalent to $g$ (hence
also ${\widetilde g}$) being locally a product of two
K{\"a}hler--Einstein metrics. {\bf Q.E.D.}

\subsection{Curvature obstructions to existence of
strictly almost K{\"a}hler structures}
The proof of Theorem \ref{th0} will be derived by showing the integrability
of
the almost \ka structure obtained in Corollary 1.
To do this we first observe that existence of a strictly almost K{\"a}hler
structure imposes several
non-trivial relations between different $U(n)$-components of the curvature.
Because
the almost K{\"a}hler structure will take the center stage in what follows,
we drop the bar-notation from the previous sub-section
and will even forget for now that in our situation
the manifold also admits a K{\"a}hler structure.

Thus, let $(M, g, J)$ be an
almost K{\"a}hler manifold  of (real) dimension $2n$.
We start by reviewing some necessary elements of
almost K{\"a}hler geometry.

The almost complex structure $J$ gives rise
to a type decomposition of complex vectors and forms, and accordingly,
of any complex tensor field;
by convention, $J$ acts on the cotangent bundle $T^*M$ by $Ja(X) = a(-JX)$.
We thus have a decomposition of the
complexified cotangent bundle
$$T^*M\otimes {\Bbb C}= \La^{1,0}M \oplus \La^{0,1}M ,$$
and of the bundle of complex
2-forms
$$\La^2M\otimes {\Bbb C}=
\La^{1,1}M \oplus \La^{2,0}M \oplus \La^{0,2}M .$$
A similar decomposition holds
for the complex bundle $S^2M\otimes {\Bbb C}$ of symmetric 2-tensors.
When considering {\it real} sections of $\La^2M$ (resp. of $S^2 M$),
we prefer to introduce
the super-scripts $'$ and $''$ for denoting the
projections to
the real sub-bundles $\La^{1,1}_{\Bbb R} M$ (resp. $S^{1,1}_{\Bbb R} M$) of
$J$-invariant 2-forms (resp. symmetric 2-tensors)
and  to $\leftr \La^{2,0}M \rightr $ (resp. $\leftr S^{2,0}M \rightr$) of
$J$-anti-invariant ones (here and henceforth $\leftr \ \ \rightr$ stands for
the real vector bundle underlying  a given complex bundle). Thus, for any
section $\psi$ of $\La^2 M$ (resp. of $S^2 M$) we have the splitting
$\psi = \psi' + \psi'' ,$ where
$$\psi'(\cdot , \cdot) = \frac{1}{2}(\psi(\cdot , \cdot) + \psi(J\cdot,
J\cdot)) \; \mbox{ and } \;
\psi''(\cdot , \cdot) = \frac{1}{2}(\psi(\cdot , \cdot) - \psi(J\cdot,
J\cdot)) .$$
Note that $\La^{1,1}_{\Bbb R}M$ can  be identified with
$S^{1,1}_{\Bbb R}M$
via the complex structure $J$:
for any $\alpha \in \La^{1,1}_{\Bbb R} M$,
$$ A = (J\circ \alpha) : = \alpha(J\cdot, \cdot) $$
is the corresponding section of $S^{1,1}_{\Bbb R}M$.

The real bundle $\leftr \La^{2,0}M \rightr$ (resp. $\leftr S^{2,0}M\rightr $)
inherits a canonical complex structure ${J}$, acting by
$$({J}\psi)(X,Y) := -\psi(JX,Y),
\ \forall \psi \in \leftr \La^{2,0}M \rightr . $$
 (We adopt a similar
definition for the action of ${J}$ on
$\leftr S^{2,0}M \rightr$).

It is well known that the fundamental form
$\Omega(\cdot,\cdot)=g(J\cdot,\cdot)$ of an almost K{\"a}hler structure
is a real harmonic 2-form of type $(1,1)$,
i.e. satisfies:
$$\Om(J\cdot, J\cdot)= \Om(\cdot,\cdot) \; \; ,  {\rm d} \Om = 0 \mbox{ and
} \delta \Om = 0 ,$$
where ${\rm d}$ and $\delta$ are the
differential and co-differential operators acting on forms. Moreover,
if $\na$ is the Levi-Civita connection of $g$, then $\na \Om$
(which is identified with the Nijenhuis tensor of $J$) is a section
of the  real vector bundle $\leftr \La^{1,0}M \otimes \La^{2,0}M \rightr $.

We first derive several consequences from
the classical Weitzenb{\"o}ck formula for a 2-form $\psi$:
\begin{eqnarray} \label{wtz2}
 \Delta \psi -\na^* \na \psi = [{\rm Ric}(\Psi \cdot,\cdot) -  {\rm
Ric}(\cdot, \Psi\cdot)] -
2 {R}(\psi) \\ \nonumber
= \frac{2(n-1)}{n(2n-1)}s \psi -2{W}(\psi) +
\frac{(n-2)}{(n-1)}[{\rm Ric}_0(\Psi\cdot,\cdot) -  {\rm Ric}_0(\cdot,
\Psi\cdot)],
\end{eqnarray}
where: $\Delta = {\rm d} \delta + \delta {\rm d}$ denotes the
 Riemannian Laplace operator
 acting on 2-forms, $\na^*$ denotes the adjoint of $\na$ with respect to $g$;
${\rm Ric}_0= {\rm Ric} -\frac{s}{2n}g$ is the traceless part of
the Ricci tensor,  $s= {\rm trace}({\rm Ric})$ is the scalar curvature,
$\Psi$ is the skew-symmetric endomorphism of $TM$ identified to $\psi$
 via the metric, and
$R$ and $W$ are respectively
the curvature tensor and the Weyl tensor, considered
as  endomorphisms of $\La^2 M$ or as sections
of $\La^2M \otimes \La^2M$, depending on the context.

Applying relation (\ref{wtz2}) to the (harmonic) fundamental form $\Om$
of the
almost K{\"a}hler structure $(g,J)$, we obtain
\begin{equation} \label{wtz2Om}
 \na^* \na \Om = 2{R} (\Om) - [{\rm Ric}(J\cdot, \cdot) -{\rm Ric}(\cdot,
J\cdot )] \; .
\end{equation}

Note that
the Ricci tensor of a K{\"a}hler structure is $J$-invariant, but this is
no longer
true for an arbitrary almost K{\"a}hler structure. It will be thus
useful to introduce
the invariant and the anti-invariant parts of the Ricci tensor with respect to
the almost complex structure $J$, ${\rm Ric}'$ and ${\rm Ric}''$,
respectively. We also put
$$\rho= J\circ {\rm Ric}'$$
to be the (1,1)-form
corresponding the the $J$-invariant part of ${\rm Ric}$, which will be called
{\it Ricci form} of $(M,g,J)$. For K{\"a}hler manifolds,
$\rho$ is clearly equal to the image of $\Om$ under the action of the
curvature $R$,
but this is not longer true for almost K{\"a}hler manifolds. In fact,
$$\rho^* = R(\Om)$$
can be considered as  a second (twisted) Ricci form of $(M,g,J)$ which is not,
in general,  $J$-invariant (see e.g. \cite{T-V}).
We will consequently denote by $(\rho^*)'$ and
$(\rho^*)''$ the corresponding 2-forms which are sections of the bundles
$\La^{1,1}_{\Bbb R} M $  and
$\leftr \La^{2,0}M \rightr $, respectively.
With these notations, formula (\ref{wtz2Om}) is a measure of the
difference of the two types of Ricci forms on an almost K{\"a}hler manifold:
\begin{equation} \label{r*r}
\rhs - \rho = \frac{1}{2}(\na^* \na \Om ). \;
\end{equation}
Taking the inner product with $\Om$ of the relation (\ref{r*r})
we obtain the difference of the two types of scalar curvatures:
\begin{equation} \label{s*s}
s^* - s = |\na \Om|^2 = \frac{1}{2} |\na J|^2 \; ,
\end{equation}
where, we recall $s= {\rm trace} ({\rm Ric})$ is the usual scalar curvature of $g$,
and
$s^* = 2\langle R(\Om),\Om \rangle$ is the so-called
{\it star-scalar curvature} of
the almost K{\"a}hler structure $(g,J)$.
Here and throughout the paper, the
inner product induced by the metric
$g$ on various tensor bundles over the manifold will be denoted by $\langle
\ \ \rangle $,
while the corresponding norm
is denoted by $| \ \ |$;
note that $\langle  \ \ \rangle $  acting on 2-forms
differs by a factor of 1/2 compared to when it acts on
corresponding tensors or endomorphisms.
In the present paper $\na \Om$ is viewed as a
$\La^2 M$-valued 1-form, while $\na J$ is considered as a section of
$(T^*M)^{\otimes 2}\otimes TM$, etc.

Formulae (\ref{r*r}) and (\ref{s*s}) can be interpreted as
``obstructions'' to the (local) existence of a strictly almost
K{\"a}hler structure $J$, compatible with a given metric $g$; see
e.g. \cite{Arm}. We derived these relations by using properties
of the  2-jet of $J$ (although eventually (\ref{s*s}) depends on
the 1-jet only), so that (\ref{r*r}) and (\ref{s*s}) can be viewed as
obstructions to the lifting of the 0-jet of $J$ to the 2-jet.

\vspace{0.2cm}
\noindent
{\it Remark 2.} In the vein of what was mentioned above one can easily
derive local non-existence results of compatible strictly almost K{\"a}hler
structures
for certain Riemannian metrics:
For example,
if we denote by
$P$ the curvature type operator acting on $\La^2M$ by
$$P(\psi) = \frac{2(n-1)}{n(2n-1)}s \psi -2{W}(\psi),$$
then, by (\ref{wtz2}) and (\ref{s*s}),
$s - s^*= \langle P(\Om), \Om \rangle = -|\na \Om|^2\le 0$.
This shows that
Riemannian metrics for which ${P}$ is semi-positive definite
do not admit even locally compatible strictly almost
K{\"a}hler structures; the latter curvature condition is equivalent
to the non-negativity of the {\it isotropic} sectional curvatures.
This criterion of non-existence applies in particular to conformally
flat manifolds of non-negative scalar curvature, or
to the symmetric spaces of {\it compact type}; see \cite{Ga,her-lamoneda}.

\vspace{0.2cm}

\noindent
In fact, there is even a more general identity than (\ref{r*r}), due
to Gray \cite{Gray},
which could also be interpreted as an obstruction
to the lifting of the 0-jet of $J$ to the 2-jet: Starting from the splitting
$$\La^2M = \La^{1,1}_{\Bbb R}M \oplus \leftr \La^{2,0}M \rightr , $$
we denote by ${\widetilde R}$ the component of the curvature operator
acting trivially on the first factor, i.e.
$${\widetilde R}_{X,Y,Z,T} = \frac{1}{4}\Big (
R_{X,Y,Z,T}- R_{JX,JY,Z,T} - R_{X,Y,JZ,JT} + R_{JX,JY,JZ,JT} \Big).$$
Thus, ${\widetilde R}$ can be viewed as a section of
the bundle ${\rm End}_{\Bbb R}(\leftr \La^{2,0}M \rightr)$, which in turn
decomposes
further as
$${\rm End}_{\Bbb R}(\leftr \La^{2,0}M \rightr) =
\Big ({\rm End}_{\Bbb R}(\leftr \La^{2,0}M \rightr)\Big)' \oplus
\Big( {\rm End}_{\Bbb R}(\leftr \La^{2,0}M \rightr)\Big)'' \; ,$$
into the sub-bundles of endomorphisms of $\leftr \La^{2,0}M \rightr$ which
commute,
respectively, anti-commute
with the action of ${J}$ on $\leftr \La^{2,0}M \rightr$. Denoting by
${\widetilde R}'$ and
${\widetilde R}''$ the corresponding components of ${\widetilde R}$,
Gray's identity is \cite{Gray}
\begin{equation} \label{gray}
 {\widetilde R}' =
- \frac{1}{4} \sum (\na_{e_j}\Om) {\otimes} (\na_{e_j}\Om) \; .
\end{equation}
As for the component ${\widetilde R}''$, from its definition we have
$${({\widetilde R}'')}_{X,Y,Z,T} =
\frac{1}{8}\Big( R_{X,Y,Z,T} - R_{JX,JY,Z,T} - R_{X,Y,JZ,JT} +
R_{JX,JY,JZ,JT}$$
$$ \hspace{3cm} + R_{X,JY,Z,JT} + R_{JX,Y,Z,JT} + R_{X,JY,JZ,T} +
R_{JX,Y,JZ,T}\Big),$$
showing that
$$({\widetilde R}'')_{Z_1,Z_2,Z_3,Z_4} =
R_{Z_1,Z_2,Z_3,Z_4}= W_{Z_1,Z_2,Z_3,Z_4} \ \ \forall Z_i \in T^{1,0}M. $$
Thus,  ${\widetilde R}''$ is actually determined by the Weyl curvature of $M$.

The next result provides a further obstruction, this time to the lift of the
3-jet of $J$ to the 4-jet (see also \cite{AA} for more detailed discussion).
\begin{prop} \label{prop2} For any almost \ka structure
$(g, J, \Om)$ the following
relation holds:
\begin{eqnarray} \label{laplace1}
\Delta(s^* - s) &=& -4 \delta \big(J \delta (J\ranti )\big) +
8 \delta \big(\langle \rho^* , \na_{\cdot} \; \Om\rangle\big) + 2|{\rm
Ric}''|^2\\ \nonumber
 & &  - 8|{\widetilde R}''|^2  - |\na^*\na \Om |^2 - |\phi|^2
+ 4\langle {\rho}, \phi \rangle - 4\langle {\rho}, \na^*\na \Om \rangle  \; ,
\end{eqnarray}
where the semi-positive (1,1)-form  $\phi$ is given by
$\phi(X,Y) = \langle \na_{JX} \Om,\na_Y \Om \rangle $;
$\delta$ denotes the co-differential with respect to $\na$, acting on
1-forms and $2$-tensors.
\end{prop}

Integrating (\ref{laplace1}) over the manifold, one obtains an integral formula identical to the one in \cite[Proposition 3.2]{sek}, up to some integration by parts. In particular we have:
\begin{cor}\label{cor3} {\rm (\cite{sek})}
For any compact almost \ka manifold  with $J$-invariant Ricci tensor
the following inequality holds:
\begin{equation} \label{intineq}
\int_M \left[ 4 \langle {\rho}, \phi \rangle  -
4 \langle {\rho}, \na^*\na \Om \rangle  - |\na^*\na \Om |^2
- |\phi|^2 \right]
dV_g \ge 0,
\end{equation}
where $dV_g= \frac{1}{n!}\Omega^{n}$ is the volume form of $g$.
\end{cor}

\noindent {\it Remark 3.} As shown by Sekigawa, the above
inequality gives an obstruction to the (global) existence of
strictly almost K{\"a}hler structures, when the metric $g$ is
Einstein with non-negative scalar curvature. Indeed, in this case
$\ranti = 0$ and $2 \langle {\rho}, \phi \rangle = \langle
{\rho}, \na^*\na \Om \rangle = \frac{s}{2n}|\na \Om|^2$, so that,
by (\ref{intineq}), $\na \Om =0$, i.e. $J$ is necessarily
K{\"a}hler. In dimension 4, other integrability results have been
derived from (\ref{laplace1}); see e.g.   \cite{Dr}, \cite{AA},
\cite{Og-Sek}.

\vspace{0.2cm}

\noindent
{\bf Proof of Proposition \ref{prop2}.} We start from the identity
$ \frac{1}{2} s^* = \langle  \Om \otimes \Om, R\rangle$, which follows from
the very definition of $s^*$.
(Note that on the bundle $\La^2 M \otimes \La^2 M$ the inner product
$\langle \ \ \rangle$ induced by $g$ differs now by a factor of 1/4
from the one induced on $(T^*M)^{\otimes 4}$.)
Applying the Laplacian to both sides of the above equality, we get
\begin{eqnarray}\label{lapl2}
 \frac{1}{2} \Delta s^* &=&
                   \Delta \langle  \Om \otimes \Om, R\rangle \\ \nonumber
                        &=&
       -2\sum \Big(\nabla_{e_i} \langle \nabla_{e_i} (\Om \otimes \Om), R
\rangle \Big)  \\ \nonumber
                  & &  - \langle \na^* \na (\Om \otimes \Om), R\rangle +
                  \langle \Om \otimes \Om , \nabla^* \nabla R \rangle.
   \end{eqnarray}
We next consider each of the terms appearing in the right-hand side of relation
(\ref{lapl2}).
Since $R$ is a {\sl symmetric} endomorphism of $\La^2 M$, the first term gives
\begin{eqnarray} \label{1st}
-2\sum \Big(\nabla_{e_i} \langle \nabla_{e_i} (\Om \otimes \Om), R \rangle
\Big) &=&
 -4 \sum \Big(\nabla_{e_i} \langle R(\Om), \na_{e_i} \Om \rangle \Big) \\
\nonumber
 &=& \ \  4  \delta \big(\langle \rho^* , \na_{\cdot} \; \Om\rangle\big) \;.
\end{eqnarray}
Using again that $R \in S^2(\La^2 M)$, the second term can be written as
\begin{equation} \label{2nd0}
-  \langle \na^* \na (\Om \otimes \Om), R\rangle =
2\sum \langle (\nabla_{e_i} \Om \otimes \nabla_{e_i} \Om ) , R \rangle  -
2 \langle (\na^* \na \Om) \otimes \Om, R \rangle .
\end{equation}
Since $ \sum (\nabla_{e_i} \Om \otimes \nabla_{e_i} \Om )
\in \Big ({\rm End}_{\Bbb R}(\leftr \La^{2,0}M \rightr)\Big)'$,
Gray's relation (\ref{gray}) implies
\begin{equation} \label{2nd1}
 2\sum \langle (\nabla_{e_i} \Om \otimes \nabla_{e_i} \Om ) , R \rangle  =
 2\sum \langle (\nabla_{e_i} \Om \otimes \nabla_{e_i} \Om ) ,
 {\widetilde R}' \rangle  =
 - |\phi|^2 \; .
\end{equation}
From (\ref{r*r}), we  have
$$ \rhsskew = (\na^* \na \Om)'' \; \mbox{ and } \rhssym - \rho = (\na^* \na
\Om)' \; .$$
Taking these into account, we have
\begin{equation} \label{2nd2}
  2 \langle (\na^* \na \Om) \otimes \Om, R \rangle  =
2 \langle \rho , \na^* \na \Om \rangle  + |\na^* \na \Om |^2 \; .
\end{equation}
Thus from (\ref{2nd0}), (\ref{2nd1}), (\ref{2nd2}), we obtain for the
second term
\begin{equation} \label{2nd}
- \langle \na^* \na (\Om \otimes \Om), R\rangle = - |\phi|^2 - |\na^* \na
\Om |^2  -
2 \langle \rho , \na^* \na \Om \rangle  \; .
\end{equation}
To compute the last term, we start by writing the Weitzenb{\"o}ck
formula for sections of $\La^2 M \otimes \La^2 M$,
applied to the curvature
tensor $R$ (see \cite[Proposition 4.2]{Bour}):
\begin{equation} \label{wtzR}
\na^* \na R = (d^{\na} \delta^{\na} + \delta^{\na} d^{\na}) R
     - 2 R \circ(\frac{1}{2} {\rm Ric} \ocap  g - R) +
K( \stackrel{\circ}{{\rm R}} \circ \stackrel{\circ}{{\rm R}} ) \; .
\end{equation}
The notations in the above relation follow the ones used by Bourguignon in
\cite{Bour}, namely:
\begin{enumerate}
\item[$\bullet$]
$d^{\na}$ and $\delta^{\na}$ are the differential and the co-differential
operators defined  on the bundle of $\La^2 M$-valued 2-forms using $\na$;
\item[$\bullet$]
$ \ocaptext $ denotes the Kulkarni-Nomizu product which allows us to
consider any symmetric tensor $S$ as an element $S \ocap  g$ of
$\La^2M\otimes \La^2M$;
\item[$\bullet$]
For any section $A$ of $\La^2 M \otimes \La^2 M$, $\stackrel{\circ}{{\rm
A}}$ denotes the endomorphism of
$(T^* M)^{\otimes 2}$  defined by
$$(T^*M)^{\otimes 2} \ni h_{X,Y} \longmapsto ({\stackrel{\circ} {\rm
A}}(h))_{X,Y} =
\sum A(e_i,X,e_j,Y) h(e_i,e_j);$$
\item[$\bullet$]
$K(\cdot)$ is 4 times the projection of $(T^*M)^{\otimes 4}$ onto
$\La^2 M \otimes \La^2M$.
\end{enumerate}
Note that apparently there is  difference compared to the formula
of Proposition 4.2 in \cite{Bour}, namely the coefficient ``2''
that we have in front of the second term. That is needed because
originally the formula in \cite{Bour} was written for sections of
$(T^*M)^{\otimes 2} \otimes \La^2 M$ and there is a difference of
a factor of 2 in the definition there of $\hat{A} \circ \hat{B}$
compared to the composition in $ {\rm End} (\La^2 M)$ used here.

Using (\ref{r*r}), we obtain
\begin{eqnarray} \label{4th2}
 - 2 \langle  R \circ(\frac{1}{2} {\rm Ric} \ocap  g - R) ,
\Om \otimes \Om \rangle = 2 \langle  ( R - \frac{1}{2} {\rm Ric} \ocap  g )
(\Om),
R(\Om) \rangle  \\ \nonumber
     =  2 \langle \rhs , \rhs - \rho \rangle
        = \langle \rho , \na^* \na \Om \rangle + \frac{1}{2} |\na^* \na \Om
|^2 \; .
\end{eqnarray}
From the definitions of $K( \stackrel{\circ}{{\rm R}} \circ
\stackrel{\circ}{{\rm R}} )$ and ${\widetilde R}$, it follows
\begin{eqnarray} \label{4th3}
 \langle  K( \stackrel{\circ}{{\rm R}} \circ \stackrel{\circ}{{\rm R}} ) ,
\Om \otimes \Om \rangle &=&
 4 \langle {\widetilde R}, J{\widetilde R}J
 \rangle \\ \nonumber
 &=& 4 \Big( |{\widetilde R}'|^2 -
             |{\widetilde R}''|^2 \Big) =
  \frac{1}{2} |\phi|^2 - 4 |{\widetilde R}''|^2 \; ,
\end{eqnarray}
where for the last step we used Gray's relation (\ref{gray}).

Finally, we express the term
$\langle (d^{\na} \delta^{\na} + \delta^{\na} d^{\na}) R , \Om \otimes \Om
\rangle $.
Because of the second Bianchi identity, $d^{\na} R = 0$ and
$$  (\delta^{\na} R)_X (Y, W) = (\na_Y {\rm Ric})(X,W) - (\na_W {\rm
Ric})(X,Y) \; .$$
Further, a short computation gives
$$   (d^{\na} \delta^{\na} R)_{X,Y,W,Z} = K(A)_{X,Y,Z,W} \; ,$$
where $A_{X,Y,Z,W} = (\na^2_{X,W} {\rm Ric})(Y,Z)$.
Using this and the product rule for the differential several times, we
eventually obtain
\begin{eqnarray} \label{4th1}
 \langle  d^{\na} \delta^{\na} R , \Om \otimes \Om \rangle  &=&
\frac{1}{2} \Delta s - 2\delta (J \delta (J \ranti)) + |\ranti|^2 \\
\nonumber
 & & + 2 \langle \rho,\phi \rangle  - \langle \rho , \na^* \na \Om \rangle \; .
\end{eqnarray}
Summing up (\ref{4th1}), (\ref{4th2}), (\ref{4th3}) and using (\ref{wtzR}),
we finally get
the last term of the right hand-side of (\ref{lapl2}),
\begin{eqnarray} \label{4th}
 \langle  \nabla^* \nabla R ,  \Om \otimes \Om\rangle  &=& \frac{1}{2} \Delta s
 - 2\delta (J \delta (J \ranti))  \\ \nonumber
 & &  + |\ranti|^2 - 4 |{\widetilde R}''|^2
+ 2 \langle \rho,\phi \rangle  + \frac{1}{2} |\na^* \na \Om |^2 +
\frac{1}{2} |\phi|^2 \; .
\end{eqnarray}
Using (\ref{1st}), (\ref{2nd}), (\ref{4th}) back in relation
(\ref{lapl2}), we get the formula (\ref{laplace1}) claimed in the
statement. {\bf Q.E.D.}

\vspace{0.2cm} \noindent {\bf Proof of Theorem \ref{th0}.} We now
turn back to the notation used in Section 2.1. Thus, $(g,J,\Om)$
denotes the K{\"a}hler structure, with Ricci tensor having two
non-negative distinct constant eigenvalues $0\le \la < \mu$,
while $(g,{\bar J}, {\overline \Om})$ is the almost K{\"a}hler
structure constructed by Lemma \ref{lem1}; we shall also use the
(1,1)-forms $\alpha$ and $\beta$ introduced in Section 2.1, so
that we have
\begin{equation}\label{rhobar}
\Om = \alpha + \beta; \  {\overline \Om}= \alpha - \beta ; \
\rho = \la\alpha + \mu \beta; \  {\bar \rho}= \la \alpha - \mu\rho,
\end{equation}
where
$\rho$ and ${\bar \rho}$ are the Ricci forms of $(g,J)$ and $(g,{\bar J})$,
respectively.

For proving Theorem \ref{th0} it is enough to show that
${\bar J}$ is integrable (see Lemma \ref{lem1}), or equivalently, that
 $\na{\overline \Om}=0$.
The latter will be derived
from the integral  inequality stated in Corollary \ref{cor3} (see Remark 3).

Let ${\bar \phi}(X,Y)=\langle \na_{{\bar J} X}{\overline \Om},
\na_{Y}{\overline \Om} \rangle$ be the semi-positive definite (1,1)-form
with respect to ${\bar J}$, defined in  Proposition \ref{prop2}.
By (\ref{rhobar}) and using
the semi-positivity of the (1,1)-forms $\alpha$ and  ${\bar \phi}$, we get
\begin{eqnarray}\label{simpler1}
\langle {\bar \rho}, {\bar \phi} \rangle - \langle {\bar \rho}, \na^*\na
{\overline \Om}\rangle &=&
(\la -\mu)\langle \alpha, {\bar \phi}\rangle + (\mu -\la)\langle \alpha ,
\na^*\na {\overline \Om} \rangle
\\ \nonumber
 & & + \mu \langle {\overline \Om}, {
\bar \phi} \rangle - \mu
\langle {\overline \Om}, \na^*\na {\overline \Om} \rangle \\ \nonumber
 & = & (\la -\mu)\langle \alpha, {\bar \phi}\rangle +
 (\mu -\la)\langle \alpha , \na^*\na {\overline \Om} \rangle -
\frac{\mu}{2}|\na {\overline \Om}|^2
\\ \nonumber
 & \le & (\mu - \la) \langle \alpha, \na^*\na {\overline \Om} \rangle -
\frac{\mu}{2}|\na {\overline \Om}|^2.
\end{eqnarray}
Since $ \langle \alpha , \na {\overline \Om} \rangle = 0$ (because
$\alpha$ and $\na {\overline \Om}$
are of type $(1,1)$ and $(2,0)+(0,2)$, respectively), we have
$$\langle \alpha ,  \na^*\na {\overline \Om} \rangle = \langle \na \alpha ,
\na {\overline \Om} \rangle =
\frac{1}{2}|\na {\overline \Om}|^2 \; , $$
where in the last step we used that $\alpha = \frac{1}{2} (\Om + {\overline
\Om})$ and $ \Om$ is parallel.
Substituting into the inequality (\ref{simpler1}), we obtain
$$\langle {\bar \rho}, {\bar \phi} \rangle - \langle {\bar \rho}, \na^*\na
{\overline \Om}\rangle \le -\frac{\la}{2}|\na {\overline
\Om}|^2.$$ Since by assumption $\la \ge 0$, the latter inequality
shows that $\langle {\bar \rho}, {\bar \phi} \rangle - \langle
{\bar \rho}, \na^*\na {\overline \Om}\rangle$ is an everywhere
non-positive function and Corollary \ref{cor3} then implies that
$\na^*\na {\overline \Om} =0 $; after multiplying by ${\overline
\Om}$ we reach $\na {\overline \Om} =0$. {\bf Q.E.D.}

\vspace{0.2cm}

\noindent {\it Remark 4.} The inequality of Corollary 2 can
actually be used for {\it any} of the almost K\"ahler structures
$(g_t, {\bar J})$ constructed in Section 2.1 and after some
computation it is not difficult to determine the $t$ dependence
of each of the terms. One would then hope to obtain some
additional information taking limits when $t \rightarrow 0$ and
$t \rightarrow \infty$. Indeed
some conclusions can be drawn in the more general situation when
only the largest eigenvalue $\mu$ is assumed non-negative, but we
fell short of obtaining integrability in this case. We do obtain
that all components of the Nijenhuis tensor $N^{{\bar J}}$ of
${\bar J}$ are zero except maybe those of the form $\langle
N^{\bar J}(A,B), B' \rangle = - \langle N^{\bar J}(B,A), B'
\rangle$, for $A \in E_{\la}$, $B, B' \in E_{\mu}$. This is equivalent
to $E_{\la}$ being a totally geodesic distribution.

\section{Homogeneous K{\"a}hler manifolds}

In this section  we consider
connected simply connected
{\it homogeneous} K{\"a}hler manifolds $(M,g,J)$, meaning that
the group of all holomorphic isometries of $(M,g,J)$ acts transitively.
Note that for any such manifold the Ricci tensor,
which is completely determined by an invariant volume form,
coincides up to sign with the canonical {\it Hermitian form} of $(M,g,J)$
\cite{koszul}; in particular ${\rm Ric}$
has constant eigenvalues with respect to  $g$.

Any homogeneous K{\"a}hler manifold admits a holomorphic fibering
over a homogeneous bounded domain whose fiber, with the induced
K{\"a}hler structure,  is isomorphic to a direct product of a
flat homogeneous K{\"a}hler manifold and a simply connected
compact homogeneous K{\"a}hler manifold; cf.
\cite{G-PSh-V},\cite{Dor-Nak}. In this structure theorem an
important role is played by the Ricci tensor whose kernel
corresponds to the flat factor \cite{Nak}; thus, when the Ricci
tensor is non-negative, the manifold splits as the product of a
flat homogeneous manifold (corresponding to the kernel of ${\rm
Ric}$) and a compact homogeneous K{\"a}hler manifold (and thus
having positive Ricci form), see \cite{cheeger-gromoll}. As for
the case of compact homogeneous manifolds, Theorem \ref{th0}
implies

\begin{cor}\label{cor2} A compact irreducible homogeneous K{\"a}hler manifold
is either K{\"a}hler--Einstein, or else the Ricci tensor has
at least three distinct eigenvalues.
\end{cor}

Of course, the above corollary can be easily derived from the
classification of the compact homogeneous K{\"a}hler manifolds, see e.g.
\cite{wang}.

Considering non-trivial (K\"ahler) homogeneous  fibrations over bounded homogeneous
domains we obtain lots of examples of irreducible homogeneous
K{\"a}hler manifolds with two distinct eigenvalues $\la <0, \mu
\ge 0$ (for explicit examples in complex dimension 2 see e.g.
\cite{Shima1,Br,A-A-D}).

We are now going to provide (non-compact) simply connected {\it
irreducible} homogeneous K{\"a}hler manifolds $(M,g,J)$ with
Ricci tensor having two negative eigenvalues  $\la< \mu <0$.
According to Corollary \ref{cor1} these also provide {\it
complete} examples of Einstein strictly almost K{\"a}hler
manifolds.

Since ${\rm Ric}$ is negative definite, $(M,J)$ must be
holomorphically equivalent to a bounded homogeneous domain; see
\cite{Shima,Nak1, Nak}. It is a result of Vinberg, Gindikin and
Piatetskii-Shapiro \cite{V-G-PSh} that any such domain has a
realization as a {\it Siegel domain of type II}, i.e. a domain $D
= \{ (z,w) \in {\Bbb C}^n\times {\Bbb C}^m : {\rm Im} z - H(w,w)
\in \mathfrak{C} \}$, where $\mathfrak{C}$ is an open convex cone
(containing no lines) in ${\Bbb R}^n$ and $H:{\Bbb C}^m\times
{\Bbb C}^m \mapsto {\Bbb C}^n$ is a Hermitian map which is
$\mathfrak{C}$-positive in the sense that
$$H(w,w) \in {\overline {\mathfrak{C}}} -\{0\} \ \ \ \forall w \neq 0.$$
If $m=0$, then $D= {\Bbb R}^n + i \mathfrak{C}$, i.e. we obtain a {\it
Siegel domain of type I} (called also a {\it tube domain}).

Any Siegel domain  admits a simply transitive action of a
solvable subgroup $S$ of {\it affine transformations} of $D$, so
that, fixing a point $p \in D$, we can identify $S$ with $D$.
Then, the complex structure $J$ and the Bergman metric ${\tilde
g}$ on $D$  pull back to define a canonical left-invariant
K{\"a}hler--Einstein structure $({\tilde g}, J)$ on $S$. Our
purpose is to show the existence of other left-invariant
K{\"a}hler metrics on $S$, whose Ricci tensor has two
distinct eigenvalues. The construction is purely algebraic and
relies on the theory of {\it normal $j$-algebras}, see {\it e.g.}
\cite{G-PSh-V}, \cite{V-G-PSh}.

The (real) Lie algebra $\mathfrak{s}$ of $S$ is  equipped with  a
scalar product $\langle \ \  \rangle$ (coming from ${\tilde g}$)
and with a vector space endomorphism $j$ (coming from $J$). This
shows that $\mathfrak{s}$ has the structure of a {\it normal
j-algebra}, meaning that there exists a 1-form $\om$ with
$\om([jX,Y]) = \langle X, Y \rangle$ for any $X,Y \in
\mathfrak{s}$, cf. \cite{PSh}. (Abstractly, a normal $j$-algebra
$(\mathfrak{s},j,\om)$ is a real Lie algebra $\mathfrak{s}$
endowed with an endomorphism $j$ and a 1-form $\om$ satisfying
certain compatibility relations; corresponding to any such
algebra there exists a bounded homogeneous domain defining
$\mathfrak{s}$ as above.) A result of Piatetskii-Shapiro
\cite{PSh} (see also \cite{d'Atri}) describes the structure of
the root spaces of $\mathfrak{s}$: Letting $\mathfrak{n} =
[\mathfrak{s},\mathfrak{s}]$, the orthogonal complement
$\mathfrak{a}$ of $\mathfrak{n}$ in $(\mathfrak{s}, \langle \
\rangle )$ is a commutative sub-algebra and $\mathfrak{n}$ can be
represented as the direct $\langle \ \rangle$-orthogonal sum of
root spaces ${\mathfrak{n}}_{\epsilon} = \{ X \in  {\mathfrak{n}}
: [A,X]=\epsilon(A)X, \ \forall A \in \mathfrak{a} \}$; if we
denote by ${\epsilon}_1, ..., {\epsilon}_r$ the (non-zero) roots
whose root spaces are mapped into $\mathfrak{a}$ by $j$, then $r=
{\rm dim} \ \mathfrak{a}={\rm rank} \ \mathfrak{s}$, and with
proper labeling, all other roots are of the form
$\frac{1}{2}\epsilon_k, \frac{1}{2}(\epsilon_\ell \pm \epsilon_s);
1\le k \le r; 1\le \ell < s \le r$ (although some of them can be
zero). It follows that ${\rm dim}_{\Bbb R}
{\mathfrak{n}}_{\epsilon_k}=1$ and we then denote by $X_k$ a
$\langle \ \rangle$-unitary generator of
${\mathfrak{n}}_{{\epsilon}_k}$; using  the fact that
$[{\mathfrak{n}}_{\epsilon'},  {\mathfrak{n}}_{\epsilon''}]
\subset {\mathfrak{n}}_{\epsilon' + \epsilon''}$ together with
the orthogonality of the root spaces, one easily checks

\begin{Lemma}\label{criterion}
If ${\mathfrak{n}}_{\frac{1}{2}\epsilon_r}=0$, then the 2-form
$\alpha$ which is the $\langle \ \ \rangle$-dual of $X_r\wedge jX_r$,
is $j$-invariant and closed, i.e. for any $X,Y,Z \in \mathfrak{g}$
$$\alpha(jX,jY)= \alpha(X,Y)$$
$$\alpha([X,Y],Z) + \alpha([Z,X],Y) + \alpha([Y,Z],X)=0.$$
\end{Lemma}
Under the hypothesis of Lemma \ref{criterion}
we may define on $S$ (by translations) a left invariant
closed
(1,1)-form $\alpha $ and a  $J$-invariant
distribution ${\mathcal R} = {\rm span}(X_r,jX_r)$. It then
follows that
the (left-invariant)  almost complex structure ${\bar J}$ on $S$,
defined by ${\bar J}|_{{\mathcal R}^{\perp}} = J|_{{\mathcal R}^{\perp}}$ and
${\bar J}|_{\mathcal R} = -J|_{\mathcal R}$, where ${\mathcal R}^{\perp}$
denotes the
${\tilde g}$-orthogonal
complement of ${\mathcal R}$, has {\it closed} fundamental form (equal to
${\widetilde \Om} - 2\alpha$, where ${\widetilde \Om}$ is the K{\"a}hler form
of $\tilde{g}$), i.e. $(M,{\tilde g}, {\bar J})$ is an Einstein almost
K{\"a}hler  manifold.
Note that ${\bar J}$ is not integrable, provided that
$\mathfrak{s}$ (equivalently $D$) is irreducible (see Lemma \ref{lem1});
by virtue of Corollary \ref{cor1}, we have also
a family of left invariant K{\"a}hler metrics on $S$, whose Ricci tensor
has two distinct negative eigenvalues.

It only remains the question of existence of irreducible
Siegel domains satisfying
the hypothesis of Lemma \ref{criterion}. Actually, non-symmetric examples
can be found in each complex dimension greater than four, see
\cite[pp. 63-64]{PSh} and \cite[pp. 411-412]{d'Atri}.
But even more  interestingly,
they also exist amongst the classical (Hermitian) symmetric domains.

\begin{cor}\label{examples} Every irreducible Hermitian symmetric space
of non-compact type which admits a realization as a tube domain carries a
strictly almost K{\"a}hler structure commuting with the standard
K{\"a}hler structure.
\end{cor}

\vspace{0.2cm} \noindent {\it Proof.} Suppose that $(M,g,J)$ is a
Hermitian symmetric space of the non-compact type and denote by
$D$ one (of the many possible) realization of $(M,J)$ as a Siegel
domain (called in this case {\it symmetric} domain). It is well
known that $g$ is now the Bergman metric of $D$ and  the real
dimensions of the root spaces
${\mathfrak{n}}_{\frac{1}{2}\epsilon_k}$ are all equal to
$\frac{2m}{r}$ (see e.g. \cite{d'Atri-Miatello, satake}); in
particular, the condition
${\mathfrak{n}}_{\frac{1}{2}\epsilon_r}=0$ means that $(M,g,J)$
admits a realization as a tube domain (see above). The complete
list of the Hermitian symmetric spaces admitting a tube domain
realization can be found for example in \cite[pp.
114-118]{satake}. {\bf Q.E.D.}

\vspace{0.2cm}
\noindent
{\it Example 1.}
An explicit example of a symmetric strictly almost K{\"a}hler
manifold (of real dimension $2n, n\ge 3$) is given by ${\rm
SO}(2,n)/({\rm SO}(2)\times {\rm SO}(n)), \ n\ge 3.$ (It is well known
that this space  admits a realization as a tube domain).

\section{The two dimensional case -- Proof of Theorem \ref{th1}}

From now on we assume that $(M,g,J)$ is a K{\"a}hler surface with
two distinct constant eigenvalues $\la < \mu$. We note that in
real four dimensions, an almost complex structure ${\bar J}$
which commutes with $J$ (and differs from $\pm J$) induces the
orientation opposite to that of $(M,J)$; if we denote by
${\overline M}$ the smooth manifold $M$, but endowed with the opposite
orientation, then, by Lemma \ref{lem1}, ${\overline M}$ carries a
symplectic structure ${\overline \Om}$. As a matter of fact, a
way of rephrasing Lemma \ref{lem1} is to say that ${\overline
\Om}$ defines an {\it indefinite} K{\"a}hler metric on $(M,J)$,
meaning that ${\overline \Om}$ is a non-degenerate, closed (1,1)-form which induces the orientation opposite to the one of
$(M,J)$.

By Theorem \ref{th0} we can further assume $\la <0$;
in this case, using the same deformation trick as in the proof of
Corollary \ref{cor1},  one can deform our  metric  to one whose Ricci
tensor has constant eigenvalues $\la< \mu$ with $\la + \mu <0$.

The next lemma is a consequence of results of \cite{Kots} and
\cite{Pe} and relies essentially on the Kodaira classification of
compact complex surfaces (see e.g. \cite{BPV}), combined with the
non-triviality of certain Seiberg--Witten invariants on the
symplectic manifold ${\overline M}$, cf. \cite{taubes}.

\begin{Lemma}\label{petean} Let $(M,g,J)$ be a
compact K{\"a}hler surface whose Ricci tensor has two distinct
constant eigenvalues $\la < \mu $ with $\la + \mu < 0$.  Then the
manifold is one of the following:
\begin{enumerate}
\item[{\rm (i)}] a minimal ruled surface which is the
projectivization of a holomorphic rank 2 vector bundle over a
curve of genus at least 2 and $\la<0<\mu$, or
\item[{\rm (ii)}] a minimal properly elliptic surface which is an elliptic
fibration over a curve of genus at least 2 with no fibers of singular
reduction and
$\la<\mu = 0$, or
\item[{\rm (iii)}] a minimal surface of general type with ample
canonical bundle and with even and non-negative signature and $\la<\mu <0$.
\end{enumerate}
\end{Lemma}
\noindent {\it Proof.} We first recall that
for any K{\"a}hler surface $(M,g,J,\Om)$
of negative scalar curvature
we have $H^0({M}, {K}^{\otimes - m})=0$, where
$K$ denotes the canonical bundle of $({M},J)$ (cf. e.g. \cite{Yau1}).

J. Petean \cite{Pe} classified the compact complex surfaces
possibly admitting indefinite K{\"a}hler metrics; those of
K{\"a}hler type  which satisfy $H^0({M}, {K}^{\otimes - m})=0$
could be either the surfaces  described in (i) \& (ii), or
minimal surfaces of general type with even and non-negative
signature. This is also a consequence of more general
results of D. Kotschick \cite{Kots}.

As for the signs of the eigenvalues $\la,\mu$ we recall that
$\frac{1}{2\pi}\rho$ represents $c_1$ in $H^2(M,{\Bbb R})$ so that
\begin{equation}\label{c_1^2}
c_1^2(M) = \frac{1}{4\pi^2}\int_M \la \mu \ dV_g,
\end{equation}
where $dV_g = \frac{1}{2}\Om\wedge\Om$ is the volume form of $(g,J,\Om)$.
It is well known \cite{BPV} that the complex surfaces described
in (i) satisfy $c_1^2 <0$, for
those in (ii) we have $c_1^2=0$, while a
minimal surface of general type satisfies $c_1^2(M)>0$; since $\la + \mu <0$
we get the needed conclusions.

We complete the proof by observing that the minimal surfaces of
general type appearing in (iii) have ample canonical bundle since
$\rho$ is a negative definite representative of $c_1$. {\bf
Q.E.D.} \vspace{0.2cm}

\vspace{0.2cm}
\noindent
{\it Remark 5.}  Some minimal complex surfaces of general
type with positive signature (which therefore admit
no locally homogeneous structure) do admit indefinite
K{\"a}hler metrics (cf. \cite{kodaira}, \cite{atiyah}).
Unfortunately, we do not know if they admit K{\"a}hler metrics
with constant eigenvalues of the Ricci tensor.

\vspace{0.2cm}

Using a result of \cite{LeBrun} and arguments from \cite{A-D-K}
and \cite{A-A-D}, we prove the integrability of ${\bar J}$
provided that one of the constant eigenvalues is non-negative.
\begin{Lemma} \label{elliptic} Every  compact K{\"a}hler surface $(M,g,J)$
whose Ricci tensor has two distinct constant eigenvalues $\la < \mu$
with $\mu\ge 0$ is locally symmetric.
\end{Lemma}
\noindent {\it Proof.} According to Theorem \ref{th0} and Lemma
\ref{lem2}, we may assume  that the scalar curvature $s=2(\la +
\mu)$ is negative. Thus, we have to consider the two possibilities
for $(M,J)$ listed in Lemma \ref{petean}, (i)\&(ii).

(i) If $(M,J)$ is a minimal
ruled  surface as in Lemma \ref{petean}(i), then  by a result of
LeBrun \cite{LeBrun} $(M,g,J)$
is locally Hermitian symmetric. Alternatively, using (\ref{gt}),
one can deform our K{\"a}hler metric $g$ to one whose
Ricci tensor has eigenvalues $(-1,+1)$ (see Lemma \ref{lem2}), i.e.
to a scalar-flat K{\"a}hler metric;
the later is locally symmetric \cite{deBar-Burns}, and so is then $g$.

(ii) Let $(M,J)$ be a minimal properly elliptic surface. By Lemma
\ref{petean}, we have $\la < \mu =0$. Adopting the notations used
in Section 2, $E_{\la}$ and $E_{\mu}$ denote the corresponding
eigenspaces of  ${\rm Ric}$ and so on. We thus have that ${\rm
Ric}$ is semi-negative definite and ${\rm Ker}({\rm Ric})=
E_{\mu}$. From Lemma \ref{petean} (ii) we also know that $(M,J)$
is an elliptic fibration over an irrational complex curve with no
fibers of singular reduction. As observed for example in
\cite{A-D-K}, after replacing $(M,J)$ with a finite cover if
necessary, we may assume that $(M,J)$ carries a non-vanishing
holomorphic vector field $X$ (acting by translations on the
fibers). Now the Bochner--Lichnerowicz formula (cf. {\rm e.g.}
\cite{besse}) shows that $X$ belongs to $E_{\mu}$ and that $X$ is
parallel with respect to the Levi-Civita connection of $g$; it
follows that $E_{\mu}$ is parallel and so is its orthogonal
complement $E_{\la}$. Since $J$ is parallel by assumption, the
almost complex structure ${\bar J}$ defined in Lemma \ref{lem1}
must be parallel too. Thus $(g,{\bar J})$ is \ka  and  $(M,g)$ is
locally symmetric, as observed in Lemma \ref{lem1}. {\bf Q.E.D.}

\vspace{0.2cm}

Let us now consider the case (iii)
of Lemma \ref{petean}.
\begin{Lemma}\label{signature} A complex surface $(M,J)$  carries a \ka
metric $g$ whose Ricci tensor has two distinct, negative, constant
eigenvalues if and only if there exists a K{\"a}hler--Einstein
metric ${\widetilde g}$ on $(M,J)$ of negative scalar curvature,
which admits a compatible almost K{\"a}hler structure ${\bar J}$
commuting with $J$;   moreover, ${\bar J}$ is integrable
precisely when $g$ {\rm (}and ${\widetilde g}${\rm )} is locally
symmetric. If the manifold is assumed compact, ${\bar J}$ is
integrable if and only if the signature of $M$ is zero.
\end{Lemma}
\noindent {\it Proof.} The first part of the lemma follows from
Corollary \ref{cor1}; by Lemma \ref{lem1} we conclude that the
integrability of ${\bar J}$ is equivalent to $(M, g)$ (hence also
$(M, {\widetilde g})$) being locally symmetric, i.e. locally a
product of two curves of distinct constant Gauss curvatures. In
the compact case, the signature of any such complex surface  is
zero. To complete the proof it remains to show that conversely, if
the signature of $M$ is zero, then ${\bar J}$ is necessarily
integrable.

Suppose for contradiction that the signature of $M$ is zero and that
${\bar J}$ is not integrable. It then follows by Lemma \ref{lem1},
Remark 1 and \cite[Th.~2]{Dr} that
$c_1({\overline M})\cdot [{\overline \Om}] <0$.
By a result of  Taubes \cite{taubes2}
we conclude that
the Seiberg--Witten invariant (for the appropriate chamber) of the
symplectic manifold
$({\overline M}, {\overline \Om})$ is non-zero so that, according to a result
of Leung \cite{leung} (see also \cite{Kots}), $(M,J)$ is uniformized by a
polydisk. From the uniqueness of the K{\"a}hler--Einstein metric
\cite{Yau} we conclude that
${\widetilde g}$ is locally
a product metric and therefore changing the orientation on one of the factors
defines  a (locally product) K{\"a}hler structure $({\widetilde
g},{\widetilde J})$
on ${\overline M}$;
it then follows for example by
\cite{A-D} that ${\bar J}$ coincides (up to sign) with ${\widetilde J}$,
 i.e. ${\bar J}$ is integrable, a contradiction.
{\bf Q.E.D.}

\vspace{0.2cm}

\noindent
{\it Remark 6.} Claude LeBrun suggested to us an alternative way to see
the connection with the Goldberg conjecture, as follows:
Let $(g,J)$ be a K{\"a}hler structure with  Ricci tensor ${\rm Ric}^g$
having two
distinct negative eigenvalues. The metric $g$ has constant {\it central}
curvature, so considering the  K{\"a}hler metric
${\widetilde g} = - {\rm Ric}^g $, this has the same Ricci form as the metric
$g$ (see  \cite{Masch}). Thus $({\widetilde g}, J)$ is a
K{\"a}hler--Einstein metric. Now it is not hard to see that the
K{\"a}hler form of $g$  is a harmonic form of constant length with
respect to ${\widetilde g}$; so will then be its {\it anti-self-dual}
part which gives rise to a negative almost K{\"a}hler (Einstein)
structure $({\widetilde g}, {\bar J})$.

\vspace{0.2cm}
\noindent
{\bf Proof of Theorem \ref{th1}.}
Theorem \ref{th1} now follows from Theorem \ref{th0} and
Lemmas \ref{elliptic} and \ref{signature}.

\vspace{0.4cm} \noindent {\bf Acknowledgments.} We are very
grateful to Claude LeBrun for suggesting the relation between the
constant eigenvalues problem and the Goldberg conjecture, to
Gideon Maschler for explaining his results on central K{\"a}hler
metrics, and to Dimitrii Alekseevsky for many useful and
instructive discussions on bounded homogeneous domains. We would
also like to thank the referee for useful observations and to
several other people with whom we had stimulating discussions --
J. Armstrong, R. Bielawski, Ph. Delanoe, P. Gauduchon and C.
Simpson. The first-named author would like to thank the
hospitality of Max-Planck-Institut in Bonn and IH{\'E}S during
the preparation of this work.


\begin{thebibliography}{99}

\bibitem{AA} V. Apostolov, J. Armstrong, {\it Symplectic 4-manifolds
with Hermitian Weyl tensor}, Trans. Amer. Math. Soc., 352 (2000), 4501--4513.


\bibitem{A-A-D} V. Apostolov, J. Armstrong and  T. Draghici, {\it Local
rigidity of certain classes of almost K{\"a}hler 4-manifolds}, Preprint 1999,
available at math.DG/9911197.

\bibitem{A-D} V. Apostolov and T. Draghici,
{\it Hermitian conformal classes and almost K{\"a}hler structures on four
 manifolds},  Diff. Geom. Appl.,  11 (1999), 179--195.

\bibitem{A-D-K} V. Apostolov, T. Draghici and D. Kotschick, {\it An
integrability theorem for almost K{\"a}hler 4-manifolds},
C. R. Acad. Sci., Paris, Ser. I, t. 329 (1999), 413--418.


\bibitem{Arm} J. Armstrong, {\it An Ansatz for Almost-K{\"a}hler, Einstein
4-manifolds}, to appear in
J. reine angew. Math.

\bibitem{atiyah} M. Atiyah, {\it The signature of fiber bundles},
in ``Global Analysis'', Papers in Honor of K. Kodaira,
Tokyo University Press (1969).


\bibitem{deBar-Burns} P. de Bartolomeis and D. Burns, {\it Stability of
vector bundles
and extremal metrics}, Inv. Math. 92 (1988), 403--407.


\bibitem{BPV} W.Barth, C.Peters and A. Van de Ven, {Compact complex
surfaces}, Springer-Verlag, 1984.

\bibitem{besse} A. L. Besse, Einstein manifolds, Ergeb. Math. Grenzgeb.3,
Folge 10,
Springer-Verlag, Berlin, Heidelberg, New York, 1987.

\bibitem{bielawski} R. Bielawski, {\it Ricci-flat K\"ahler metrics on
canonical bundles}, preprint 2000, available at math.DG/0006144.

\bibitem{BI} D. E. Blair and S. Ianus, {\it Critical associated metrics on
symplectic manifolds},
Contemp. Math. 51 (1986), 23--29.

\bibitem{Bour} J.-P. Bourguignon, {\it Les vari{\'e}t{\'e}s de dimension 4
{\`a} signature non nulle dont la courbure est harmonique sont
d'Einstein}, Inv. Math. 63 (1981), 263-286.

\bibitem{Br} R. Bryant, private communication.

\bibitem{cheeger-gromoll} J. Cheeger and S. Gromoll, {\it The splitting
theorem for
manifolds of non-negative Ricci curvature}, J. Differential Geom., 6
(1971), 119--128.

\bibitem{d'Atri} J.E. D'Atri,
{\it Holomorphic sectional curvatures of bounded homogeneous domains and
related questions}, Trans. Amer. Math. Soc. 256 (1979), 405--413.

\bibitem{d'Atri-Miatello} J.E. D'Atri and D. Miatello, {\it A
characterization of
bounded symmetric domains by curvature},
Trans. Amer. Math. Soc., 267 (1983), 531--540.

\bibitem{DM} J. Davidov and O. Mu\u{s}karov, {\it Twistor spaces with
Hermitian Ricci tensor}, Proc. Amer. Math. Soc. 109 (1990), no. 4, 1115--1120.

\bibitem{Dor-Nak} J. Dorfmeister and K. Nakajiama,
{\it The fundamental conjecture for homogeneous K{\"a}hler manifolds},
Acta Math., 161 (1988), 23--70.

\bibitem{Dr} T. Draghici, {\it Almost \ka 4-manifolds with $J$-invariant
Ricci tensor}, Houston J. Math., 25 (1999), 133--145.

\bibitem{her-lamoneda} L. Hern{\'a}ndez-Lamoneda, {\it Curvature vs. Almost
Hermitian Structures}, Geom. Dedicata, 79 (2000), 205--218.

\bibitem{hwang-singer} A. Hwang and M. Singer, {\it A momentum construction
for circle-invariant K\"ahler metrics}, available at math.DG/9811024.

\bibitem{Ga} P. Gauduchon, {\it Complex structures on compact conformal
manifolds
of negative type}, in  ``Complex Analysis and Geometry'', (eds. V. Ancona,
E. Ballico and A. Silva) Marcel Dekker, New York-Basel-Hong Kong, 1996,
201--212.

\bibitem{G-PSh-V} S.G. Gindikin, I.I. Piatetskii-Shapiro, and E.B. Vinberg,
{\it Homogeneous K{\"a}hler manifolds}, in ``Geometry of Homogeneous
Bounded Domains'',
C.I.M.E., 1967, 1--88.

\bibitem{goldberg} S.I. Goldberg, {\it Integrability of almost K{\"a}hler
  manifolds}, Proc. Amer. Math. Soc. 21 (1969), 96-100.

\bibitem{Gray} A. Gray, {\it Curvature identities for Hermitian and almost
Hermitian manifolds}, T{\^o}hoku Math. J. 28 (1976), 601-612.


\bibitem{kodaira} K. Kodaira,
{\it A certain type irregular algebraic surfaces}, J. Analyse Math. 19 (1967),
207--215.



\bibitem{koszul} J.L. Koszul, {\it Sur la forme hermitienne canonique
des espaces homog{\`e}nes complexes}, Canad. J. Math. 7 (1955), 562--576.

\bibitem{Kots} D. Kotschick, {\it Orientations and geometrisations of compact
complex surfaces}, Bull. London Math. Soc. 29 (1997), no. 2,
145--149.

\bibitem{Kow} O. Kowalski, Generalized Symmetric Spaces, LNM 805, 1980.

\bibitem{LeBrun} C. LeBrun, {\it Polarized 4-manifolds, extremal
K{\"a}hler metrics and Seiberg-Witten theory}, Math. Res. Lett. 2, No.5,
653--662 (1995).

\bibitem{leung} N.C. Leung,
{\it Seiberg-Witten invariants and uniformizations},
Math. Ann. 306 (1996), 31--46.


\bibitem{Masch} G. Maschler, {\it Central K{\"a}hler metrics}, Preprint.

\bibitem{And} A. Moroianu, {\it On Kirchberg inequality for compact K{\"a}hler
    manifolds of even complex dimension, } Ann. Global
  Anal. Geom. {15} (1997), 235--242.

\bibitem{andrej} A. Moroianu, {\it K{\"a}hler manifolds with small eigenvalues
of the Dirac operator and a conjecture of Lichnerowicz}, Ann. Inst. Fourier
(Grenoble), 49 (1999), 1637--1659.

\bibitem{Nak} K. Nakajima, {\it Homogeneous K{\"a}hler manifolds of
non-degenerate
Ricci curvature}, J. Math. Soc. Japan, 42 (1990), 475--494.

\bibitem{Nak1} K. Nakajima, {\it Homogeneous K{\"a}hler manifolds of
non-positive
Ricci curvature}, J. Math. Kyoto Univ., 26 (1986), 547--558.

\bibitem{Og-Sek} T. Oguro, K. Sekigawa,
{\it Four-dimensional almost K{\"a}hler Einstein
and $*$-Einstein manifolds}, Geom. Dedicata, 69 (1998), 91--112.

\bibitem{Olszak} Z. Olszak,
{\it A note on almost K{\"a}hler manifolds},
Bull. Acad. Polon. Sci., XXVI (1978), 139--141.

\bibitem{Pe} J. Petean, {\it Indefinite K{\"a}hler-Einstein metrics on
compact complex surfaces}, Comm. Math. Phys. 189 (1997), 227--235.

\bibitem{PSh} I.I. Piatetskii-Shapiro, Automorphic functions and the
geometry of classical domains, English transl., Gordon and Breach, New
York, 1969.

\bibitem{satake} I. Satake, Algebraic Structures of Symmetric Domains,
{\it Iwanami Shoten Publ. and Princeton University Press}, 1980.

\bibitem{sek} K. Sekigawa,
{\it On some compact Einstein almost K{\"a}hler manifolds},
J. Math. Soc. Japan, 39 (1987), 677--684.

\bibitem{Shima} H. Shima, {\it On homogeneous K{\"a}hler manifolds of
solvable Lie groups},
J. Math. Soc. Japan, 25 (1973), 422--445.

\bibitem{Shima1} H. Shima, {\it On homogeneous K{\"a}hler manifolds
with non-degenerate canonical Hermitian form of signature $(2, 2(n-1))$},
Osaka J. Math., 10 (1973), 477--493.

%%\bibitem{Siu} Y.T. Siu, {\it Every K3 surface is K{\"a}hler}, Inv.
%%Math., 73 (1983), 139--150.

\bibitem{taubes} C. H. Taubes, {\it The Seiberg-Witten Invariants and
Symplectic Forms}, Math. Res. Lett., 1 (1994), 809--822.

\bibitem{taubes2} C. H. Taubes, {\it ${\rm SW}\Rightarrow{\rm Gr}$: from the
Seiberg-Witten equations to pseudo-holomorphic curves}, J. Amer. Math.
Soc., 9(1996), no.3, 845--918.

\bibitem{T-V} F. Tricerri and L. Vanhecke, {\it Curvature tensors on almost
Hermitian manifolds}, Trans. Amer. Math. Soc., 267 (1981), 365--398.

\bibitem{V-G-PSh} E.B. Vinberg, S.G. Gindinkin, and I.I. Piatetskii-Shapiro,
{\it Classification and canonical realization of complex bounded homogeneous
domains}, Trudy Moskov. Mat. Obsc. 12 (1963), 359--388;
[Trans. Moscow Math. Soc., 1963, 404--437].

\bibitem{wall} C. T. C. Wall, {\it Geometric structures on compact complex
analytic surfaces}, Topology 25 (1986), 119--153.

\bibitem{Yau} S.-T. Yau, {\it On the Ricci curvature of a compact
K{\"a}hler manifold and the complex Monge-Amp{\`e}re equation I,} Comm.
Pure Appl.
Math. 31, (1978) 339--411.

\bibitem{Yau1} S.-T. Yau,
{\it On the scalar curvature of compact Hermitian manifolds},
Inv. Math. 29 (1974), 213--239.

\bibitem{wang} H.C. Wang, {\it Closed manifolds with homogeneous complex
structure}, Amer. J. Math., 76 (1954), 1--32.

\end{thebibliography}
\end{document}